\documentclass[12pt]{amsart}

\usepackage{amsmath,amssymb,amsfonts,amsthm,braids,graphics,verbatim}

\usepackage[all,2cell,dvips]{xy}

\nc{\margin}[1]{\marginpar{\scriptsize #1}}

\include{diagram}

\begin{document}

\input{epsf.sty}

\title{\bf Braid groups and the co-Hopfian property}

\author{Robert W. Bell}
\author{Dan Margalit}

\address{Department of Mathematics\\ University of Utah\\ 155 S 1440 East \\ Salt Lake City, UT 84112-0090}

\thanks{Both authors are partially supported by a VIGRE postdoctoral position under NSF grant number 0091675 to the University of Utah.}

\email{rbell@math.utah.edu, margalit@math.utah.edu}



%%%
%%%
%%%

\maketitle

\begin{center}\today\end{center}

\begin{abstract}Let $B_n$ be the braid group on $n \geq 4$ strands.  We prove that $B_n$ modulo its center is co-Hopfian.  We then show that any injective endomorphism of $B_n$ is geometric in the sense that it is induced by a homeomorphism of a punctured disk.  We further prove that any injection from $B_n$ to $B_{n+1}$ is geometric. Additionally, we obtain analogous results for mapping class groups of punctured spheres.  The methods use Thurston's theory of surface homeomorphisms and build upon work of Ivanov and McCarthy.\end{abstract}

\section{Introduction}

The {\em braid group on $n$ strands} is the group of isotopy classes
of orientation preserving homeomorphisms of the $n$-times punctured
disk $\dn$ which are the identity on the boundary (see \cite{jb}): \[ \bn =
\pi_0(\homeop(\dn,\ddn)) \] $\bn$ is generated by \emph{half-twists};
each such generator $H_a$ is the isotopy class of a homeomorphism of
$\dn$ which switches two punctures along an arc $a$ (see
Section~\ref{background}).

A group is {\em co-Hopfian} if every injective endomorphism is an
isomorphism.  We see that $\bn$ is not co-Hopfian: any map which takes
each half-twist $H_a$ to the product $H_a \gar^t$ (where $\gar$
generates the center $\cent$ of $\bn$ and $t$ is fixed nonzero
integer) is injective, but not surjective.

Given a homeomorphism (or an isotopy class) $h$ of an oriented surface, let
$\ep(h)$ equal 1 if $h$ preserves the orientation of the surface, and
let $\ep(h)$ equal $-1$ otherwise.

\begin{mthm}\label{main1}
Let $n \geq 4$.  Any injective endomorphism $\inj$ of $\bn$ is induced by a
homeomorphism $h:\dn\to\dn$ in the following sense: there is a fixed
integer $t$ so that \[ \inj(H_a)
= H_{h(a)}^{\ep(h)} \ \gar^t \] for each half-twist $H_a$.\end{mthm}

For any homeomorphism $h$ of $\dn$ and integer $t$, there is an
injective endomorphism $\inj$ as in the theorem.  The map $\rho$ is
not surjective whenever $t \neq 0$ (nothing maps to $z$).

Braid groups are {\em Hopfian}---every surjective endomorphism is an
isomorphism:  Bigelow and Krammer showed that braid groups are linear
\cite{sb} \cite{dk}; and by well-known results of Mal'cev, finitely
generated linear groups are residually finite, and finitely generated residually finite groups are Hopfian.

Lin has shown that for $k < n \neq 4$, any homomorphism from $\bn$ to $B_k$ has cyclic image \cite{vl}.  We characterize injective homomorphisms of $\bn$ into $\bnp$ (see Section~\ref{background} for the definition of a Dehn twist):

\begin{mthm}\label{main2}Let $n \geq 4$.  Any injective homomorphism $\inj: \bn \to
  \bnp$ is induced by an embedding $h: \dn \to \dnp$ in the following
  sense: there are fixed integers $s$ and $t$ so that \[ \inj(H_a) =
  H_{h(a)}^{\ep(h)} \ T_{A}^s \gar^t \]
for each half-twist $H_a$,
 where $T_{A}$ denotes the Dehn
  twist about the curve $A = h(\ddn)$.\end{mthm}

To prove this theorem, we introduce the \emph{arc triple complex}
$\atdn$ (refer to Section~\ref{ad}) and prove that this complex is
connected for $n \geq 7$.  The cases where $n < 7$ require special
consideration.

Our starting point is Main Theorem~\ref{main3} (below), which says
that every injective endomorphism of $\bmodz$ is induced by a
homeomorphism of $\dn$.  To make sense of this we use the fact that $\bmodz$ is
isomorphic to $\pi_0(\homeop(\dn))$.  In general, the \emph{mapping
  class group} of a surface $S$ is defined as:
\[ \mods = \pi_0(\homeop(S)) \]

\begin{mthm}\label{main3}
Let $n \geq 4$.  The group $\bmodz$ is co-Hopfian, and any injective endomorphism
  $\inj$ of $\bmodz$ is induced by a homeomorphism $h$ of $\dn$ in the sense that
\[ \inj(H_a) = H_{h(a)}^{\ep(h)} \]
for each half-twist $H_a$.
\end{mthm}

We will deduce Main Theorem~\ref{main1} as a corollary of Main
Theorem~\ref{main3}
in Section~\ref{almostcohopf}.
The proofs of our Main Theorems~\ref{main3}
and~\ref{main2} follow the basic strategy of the following theorem of Ivanov and
McCarthy \cite{im}:

\begin{thm}\label{im}
If $S$ is a compact orientable surface which has positive genus and is
  not a torus with $0$, $1$, or $2$ boundary components, then $\mods$
  is co-Hopfian.  Moreover there are no injective homomorphisms $\mods
  \to \mcg(S')$ where $S'$ is the surface obtained from $S$ by
  removing an open disk.
\end{thm}

The reader should contrast the above theorem with Main
Theorem~\ref{main2}, where there do exist injective homomorphisms
between different braid groups.

\nc{\Z}{\mathbb{Z}}

Note that $B_3 / \cent \cong \pslz \cong \Z_2 * \Z_3$ is not co-Hopfian.  Crisp and
Paoluzzi have recently established the co-Hopfian property for $B_4 /
\cent$, by showing that it has an essentially unique CAT(0) structure
\cite{cp}.

As another corollary of Main Theorem~\ref{main3}, we obtain the analog
of the theorem of Ivanov and McCarthy for genus zero surfaces:

\begin{thm}\label{mcgthm}If $S$ is a sphere with $n \geq 5$ punctures,
  then $\mods$ is co-Hopfian, and there are no injective homomorphisms $\mods \to \mcg(S')$, where $S'$ is a sphere with $n+1$ punctures.\end{thm}

The results of this paper can be used to recover the classical
theorems that $\out(\bn) \cong \zt$ and $\out(\mcg(S)) \cong \zt$ for
$S$ a sphere with $n+1$ punctures ($n \geq 4$).  Both were originally obtained by Dyer and Grossman via algebraic methods \cite{dg}, and later proven by Ivanov using Thurston's theory \cite{ni}.  In related work, Korkmaz proved that the abstract commensurator of $\bmodz$ is isomorphic to $\aut(\bmodz)$, and that any endomorphism of $\bmodz$ with finite-index image is an automorphism ($n \geq 4$) \cite{mk}.

A key phenomenon which allows us to promote algebraic embeddings to topological ones is the fact that all algebraic braid relations between half-twists look topologically the same (Lemma~\ref{aba}).

\p{Outline.}Section~\ref{almostcohopf} gives the proof of Main
Theorem~\ref{main1}, assuming Main Theorem~\ref{main3}. In
Section~\ref{background} we give definitions and basic constructions
used in the proofs of Main Theorems \ref{main2} and~\ref{main3}.  Section~\ref{relations} explores the interplay between the algebra and topology of half-twists in $\bn$.

Section~\ref{bmodztobmodz} is the proof of Main
Theorem~\ref{main3}. There are three main steps:
\bl
\item Step 1: Any injective endomorphism $\inj$ of $\bmodz$ is almost half-twist preserving: it takes some power of any half-twist to a power of a half-twist
\item Step 2: $\inj$ is actually half-twist preserving
\item Step 3: $\inj$ is induced by a homeomorphism of $\dn$\el

In Section~\ref{bntobnp}, we prove Main Theorem~\ref{main2}
using the same strategy.
We explain how to translate the main theorems to mapping class groups of punctured spheres in Section~\ref{mcgpf}.

We then define the arc $k$-tuple complex in Section~\ref{ad} and prove that it is connected (this is used in Section~\ref{bntobnp}).  In Section~\ref{q} we ask questions related to this work.

\p{Acknowledgements.}We wish to thank Jason Behrstock, Mladen Bestvina, Tara Brendle, Martin Bridson, Ryan Budney, Fred Cohen, Benson Farb, Edward Formanek, Steve Gersten, Allen Hatcher, Nikolai Ivanov, Richard Kent, Mustafa Korkmaz, William Thurston, Domingo Toledo, and Kevin Wortman for helpful conversations and for comments on earlier drafts.  We are indebted to Joan Birman for directing us to her paper with Hugh Hilden, and for many insightful and enthusiastic correspondences.

Near the start of this work, Vladimir Lin sent us examples of injective endomorphisms of the braid group which are not surjective.  Chris Leininger was extremely generous with his time and energy in reading our paper, suggesting modifications, and working on the special cases of the second main theorem.

The second author is especially grateful to Joan Birman, Benson Farb, and Kevin Wortman for encouragement on this project.

%%%
%%%
%%%

\section{$\bn$ is almost co-Hopfian}\label{almostcohopf}

Assuming Main Theorem~\ref{main3}, we now prove Main
Theorem~\ref{main1}.

\bpf

Since $\binj(Z) < Z$ (by Lemma~\ref{z} below) and $\binj^{-1}(Z) < Z$
(as $\binj$ is injective), there is a well-defined injective homomorphism $\binjz: \bmodz \to \bmodz$.  By Main Theorem~\ref{main3},  $\binjz$ is an isomorphism induced by a homeomorphism $h$ of $\dn$:
\[ \binj(H_a) Z = \binjz(H_a Z) = H_{h(a)}^{\ep(h)} Z \]
\[ \binj(H_a) = H_{h(a)}^{\ep(h)} \gar^{t_a} \]

The exponent $t_a$ is independent of $a$, since all half-twists are
conjugate.

\epf

We will require the following theorem of de~Ker\'ekj\'art\'o, Brouwer, and Eilenberg \cite{bdk} \cite{lejb} \cite{se}:

\begin{thm}\label{eilenberg}Any finite order homeomorphism of a disk or a sphere is conjugate (via a homeomorphism) to a Euclidean isometry.\end{thm}

The Nielsen realization theorem says that any finite order element of
 $\mcg(\dn)$ can be realized by a finite order homeomorphism \cite{spk}.
 Since $z$ is the identity as an element of $\mcg(\dn)$, the following
 is a corollary to Theorem~\ref{eilenberg}:

\begin{cor}\label{periodic}Any root of a central element of $\bn$ is conjugate to a power of one of the elements $\de$ or $\ga$ of Figure~\ref{gdpic}.\end{cor}

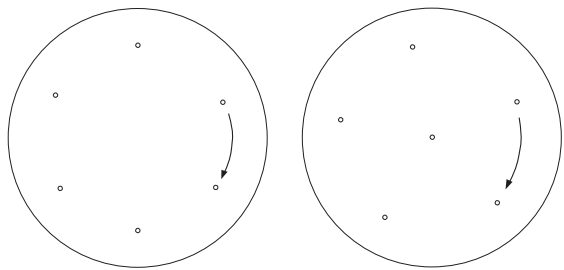
\pics{gd}{The roots $\de$ and $\ga$ of $\gar$.}{4}

\begin{lem}\label{z}If $\binj: \bn \to \bn$ is an injective homomorphism, then $\binj(Z) < Z$.\end{lem}

\bpf

Let $G$ be a free abelian subgroup of maximal rank containing $z$. Then $\binj(G)$ is also of maximal rank, and hence $\binj(G) \cap Z$ is nontrivial.

Since $\binj^{-1}(\cent) < \cent$, we have $\binj(\gar)^k \in \cent$ for some $k$.  By Corollary~\ref{periodic}, $\binj(\ga)$ and $\binj(\de)$ are conjugate to powers of $\de$ and $\ga$.  Also observe that $\de^k$ only fixes a puncture of $\dn$ if $k$ is a multiple of $n$, whereas every power of $\ga$ fixes at least one puncture of $\dn$.

If $\binj(\de)$ is conjugate to a power of $\de$, then $\binj(\gar)=\binj(\de^n)$ is central.  Similarly, if $\binj(\ga)$ is conjugate to a power of $\ga$, then $\binj(\gar)=\binj(\ga^{n-1})$ is central.  Thus, we can assume that $\binj(\de)$ is conjugate to a power of $\ga$, and vice versa.

In this case, $\binj(\gar)$ is conjugate to both a power of $\ga$ and a power of $\de$.  But by considering fixed punctures, the only conjugate powers of $\ga$ and $\de$ are central.  Therefore, $\binj(\gar) \in \cent$.

\epf

\p{Remark.} The braid groups exhibit a general obstruction to the co-Hopfian property: if $G$ is a group with a homomorphism $L : G \to \z$ and an infinite order central element $z$ with $L(z) \notin \set{0,-1,-2}$, then the endomorphism given by: \[ g \mapsto g z^{L(g)} \] is injective but not surjective.  Finite-type Artin groups are examples, where $L$ is taken to be the usual {\em length homomorphism}.

%%%
%%%
%%%

\section{Background}\label{background}

In this section we introduce ideas from mapping class groups (Thurston theory) and explain their connection to braid groups.  Throughout the paper, we use functional notation for words in $\bn$.

\p{Curves and arcs.} The {\em interior} of a simple closed curve in $\dn$ is the component of its complement which does not contain the boundary circle of $\dn$.  A simple closed curve in $\dn$ is {\em nontrivial} if it contains more than one puncture, but not all $n$ punctures, in its interior.  A {\em 2-curve} is a simple closed curve with exactly two punctures in its interior.
By an {\em arc} in $\dn$, we always mean a simple arc connecting two
distinct punctures (its {\em ends}). 

When convenient, we confuse curves and arcs with their isotopy classes.

There is a bijection of isotopy classes: \[ \text{\{2-curves\}}
\longleftrightarrow
\text{\{arcs\}} \]

We use $\gin(a,b)$ to denote the {\em geometric intersection number} of curves $a$ and $b$.

\p{Dehn twists.} The \emph{Dehn twist} about a curve $a$ is the
mapping class (isotopy class of homeomorphisms) $T_a$ whose support is an annular neighborhood of $a$ and whose action on the annular neighborhood is described by Figure~\ref{dtpic}.  The center of $\bn$ is generated by $z$, the Dehn twist about $\ddn$.

\pics{dt}{Dehn twist about a curve $a$.}{4}

A {\em multitwist} is a product of powers of Dehn twists about disjoint curves.

\p{Half-twists.} By a half-twist $H_a$ along an arc $a$, we mean the mapping class with support the interior of a 2-curve, as described in Figure~\ref{htpic}.

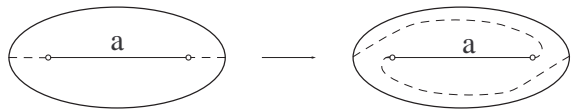
\pics{ht}{Half-twist about an arc $a$.}{2}

In light of the bijection between 2-curves and arcs, we may refer to a half-twist with respect to either an arc or a 2-curve.  Note that $H_a^2 = T_a$ for a 2-curve $a$.

\p{Adjacency.} Two disjoint arcs in $\dn$ are said to be {\em adjacent} if they share exactly one end.  This is equivalent to the condition that the corresponding 2-curves have geometric intersection number 2 (we also call such 2-curves adjacent).

\p{Conjugation.} If $f \in \mods$, then $f H_a^j f^{-1} =
H_{f(a)}^j$.  It follows that:

\begin{fact}\label{fix}For $j \neq 0$, $[f,H_a^j]=1 \iff f(a)=a$.\end{fact}

\p{Pseudo-Anosov.} By Thurston's classification, a mapping class $f$
is \emph{pseudo-Anosov} if and only if $f^k(a) \neq a$ for every
simple closed curve $a$ and any nonzero integer $k$.  This is our
working  definition of pseudo-Anosov.   By work of Ivanov, the
centralizer of a pseudo-Anosov mapping class is virtually cyclic \cite{nvi}.

\p{Reductions.} If an element $f$ of $\mods$ fixes a collection $\coll$ of disjoint curves in $S$, then there is a representative homeomorphism for $f$ which fixes a set of representatives for $\coll$.  This gives rise to a well-defined element $f_\coll$ of $\mcg(S_\coll)$, called the {\em reduction of $f$ along $\coll$}, where $S_\coll$ is the surface obtained by cutting $S$ along the representatives for $\coll$.  Note that any twist about a curve of $\coll$ becomes trivial in the reduction.

\p{Pure mapping classes.} An element $f$ of $\mods$ is called {\em pure} if it has a reduction $f_\coll$ which induces the trivial permutation on the components of $S_\coll$, acts as the identity on the boundary (fixing punctures as well), and
 restricts to either the identity or a pseudo-Anosov map on each such component (see \cite{nvi}).  In the case of a punctured disk, the set of pure mapping classes coincides with the classical {\em pure braid group} $\pbn$ (see e.g.~\cite{jb}):

\begin{thm}[Irmak-Ivanov-McCarthy \cite{iim}]\label{purity}The elements of $\pbn/\cent$ are pure in the above sense.\end{thm}

\p{Canonical reduction systems.} A curve $a$ is in the {\em canonical reduction system} for a pure mapping class $f$ if $f(a)=a$, and $f(b) \neq b$ whenever $\gin(a,b) > 0$.  This notion is due to Birman, Lubotzky, and McCarthy \cite{blm}.

We have a generalization of Fact~\ref{fix} which follows from the definition above:

\begin{lem}\label{commutecrs}If $[f,g]=1$, and $\coll$ is the canonical reduction system for $g$, then $f(\coll) = \coll$.\end{lem}

Finally, we give a modification of Lemma 10.2 of the paper of Ivanov and McCarthy (\cite{im}).  This will be used in the proofs of the main theorems.  By the {\em rank} of a group $\Ga$, denoted $\rk(\Ga)$, we mean the maximal rank of a free abelian subgroup.

\begin{lem}\label{lem 10.2 I-M} Let $\inj: \Ga \to \Ga'$ be an
  injective homomorphism.  Suppose that
$\rk \Ga' = \rk \Ga + R < \infty$ for some non-negative
integer $R$.  Let $G < \Ga$ be an abelian subgroup
of maximal rank, and let $f \in G$.  Then
\[\rk Z(C_{\Gamma'}(\inj(f))) \leq \rk Z(C_{\Gamma}(f)) + R\]
\end{lem}

\bpf

Let $A = \inj(G) \cap Z(C_{\Ga'}(\inj(f)))$, and let $B = \langle \, \inj(G), Z(C_{\Ga'}(\inj(f))) \, \rangle$.

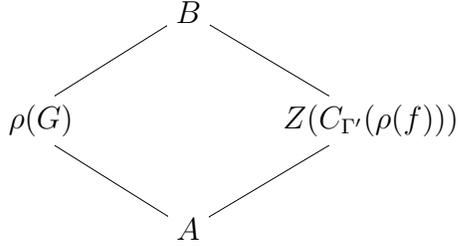
\begin{center}
\SelectTips{}{}
\hspace{.0in} \xymatrix{
 & B \ar@{-}[ld] \ar@{-}[rd] & \\
\ \ \ \ \ \ \ \inj(G) \ \ & & Z(C_{\Ga'}(\inj(f))) \\
 & A \ar@{-}[lu] \ar@{-}[ru] &
}
\end{center}

Note that all of the groups in the above diagram are abelian.
We have: \[\inj(G) \subset C_{\Ga'}(\inj(f))\]
Thus, we have: \begin{eqnarray*}
\rk \inj(G) + \rk Z(C_{\Ga'}(\inj(f))) &=&
\rk A + \rk B \\
&\leq& \rk A + \rk \Ga' \\ &=& \rk A + \rk \Ga + R\\
 &=& \rk A + \rk G + R \\
 &=& \rk A + \rk \inj(G) + R \end{eqnarray*}
Therefore, $\rk Z(C_{\Ga'}(\inj(f))) \leq \rk A + R$.

Observe that $A \subset Z(C_{\inj(\Ga)}(\inj(f)))$.  But this latter group is isomorphic to $Z(C_{\Ga}(f))$. Combined with the previous inequality, this completes the proof.

\epf

%%%
%%%
%%%

\section{Relations}\label{relations}

We use a double cover argument of Birman and Hilden to draw an analogy between braid groups and mapping class groups of higher genus surfaces.

\p{Marked points.} Let $\dpu$ be a disk with a set $\punc$ of $n$ marked points.  There is a natural isomorphism: \[ \mcg(\dn) \cong \mcg(\dpu) \] where homeomorphisms of $\dpu$ are required to fix $\punc$ as a set.

\p{Double cover.} If $n$ is odd, then there is a 2-sheeted branched cover over $\dpu$ by a surface $\dnt$ with genus $\frac{n-1}{2}$ and one boundary circle (the marked points are the branch points).  The covering transformation is an involution $\iota$ switching the two sheets.  If $n$ is even, then we take $\dnt$ to be a surface with genus $\frac{n-2}{2}$ and two boundary circles.  See Figure~\ref{tildepic}.

We note that since we only consider $\bn$ for $n \geq 4$, we have that $\dnt$ is either a torus with two boundary circles or it has genus greater than one.

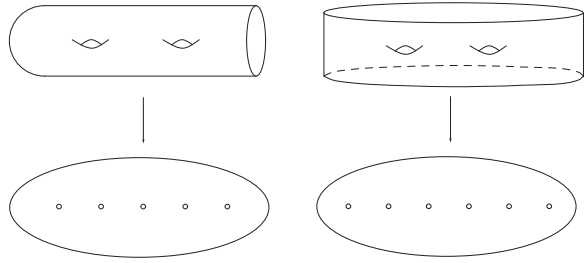
\pics{tilde}{The covering $\dnt \to \dpu$ for $n$ odd and $n$ even.}{4}

\p{Symmetric mapping class group.} Let $\hmcg(\dnt)$ denote the
centralizer of $\iota$ in $\mcg(\dnt)$.  Birman and Hilden prove the
following \cite{bh71} \cite{bh73}:

\begin{thm}\label{hyper}$\mcg(\dpu)  \cong \hmcg(\dnt) / \genby{\iota}$.\end{thm}

Birman and Hilden state this theorem for a covering of a closed surface of genus $g$ over a sphere with $2g+2$ marked points.  But their proof holds verbatim in our case (see also \cite{bh72}).

The isomorphism of Theorem~\ref{hyper} can be described explicitly on generators.  Any half-twist in $\mcg(\dpu)$ about an arc $a$ corresponds to a Dehn twist about the simple closed curve $\ta$ which is the lift of $a$ to $\dnt$.

\p{From Dehn twists...} We now use Theorem~\ref{hyper} to translate relevant facts about Dehn twists in surface mapping class groups to facts about half-twists in braid groups.  Here are the statements about Dehn twists ($j$ and $k$ are nonzero integers):

\begin{lem}\label{uniquetwist}$T_a^j=T_b^k$ $\Leftrightarrow$ $a=b$ and  $j=k$.\end{lem}

\begin{lem}\label{disjointtwist}$[T_a^j,T_b^k]=1 \Leftrightarrow \gin(a,b) = 0$.\end{lem}

\begin{lem}[Ivanov--McCarthy \cite{im}]\label{abatwist} Distinct Dehn twists $T_a$ and $T_b$ satisfy \[ T_a^jT_b^kT_a^j = T_b^kT_a^jT_b^k \] if and only if $\gin(a,b)=1$ and $j=k=\pm 1$.\end{lem}

In the next theorem, $F_2$ denotes the free group on two letters.

\begin{thm}[Ishida \cite{ai}, Hamidi-Tehrani \cite{hht}]\label{free} $\genby{T_a^j,T_b^k} \ncong F_2 \Leftrightarrow$ $\gin(a,b) = 0$ or $\gin(a,b)=1$ and $\set{j,k} \in \set{\set{1},\set{1,2},\set{1,3}}$.\end{thm}

\p{...to half-twists.} The following lemma is the desired correspondence between twist relations in surface mapping class groups and half-twist relations in braid groups.

\begin{lem}\label{word}Powers of half-twists $H_a^j$ and $H_b^k$ satisfy a relation in $\mcg(\dpu)$ if and only if the corresponding powers of twists $T_{\ta}^j$ and $T_{\tb}^k$ satisfy the same relation in $\mcg(\dnt)$.\end{lem}

\bpf

Suppose some word $W(H_a^j,H_b^k)$ in $H_a^j$ and $H_b^k$ equals the identity in $\mcg(\dpu)$.  By Theorem~\ref{hyper}, we have $W(T_{\ta}^j,T_{\tb}^k)=\iota^\ep$, so $W(T_{\ta}^j,T_{\tb}^k)^2=1$ in $\hmcg(\dnt)$ (and hence in $\mcg(\dnt)$).  Then Theorem~\ref{free} implies that $\gin(\ta,\tb) \leq 1$.  So there is a homeomorphism representing $W(T_{\ta}^j,T_{\tb}^k)$ whose support is either a pair of annuli or a torus with one boundary circle.  On the contrary, $\iota$ has no such representative (recall that $\dnt$ is ``at least'' a torus with two boundary circles), so $\ep=0$, and $W(T_{\ta}^j,T_{\tb}^k)=1$.  The other direction is trivial.

\epf

Combining Lemma~\ref{word} with Lemmas~\ref{uniquetwist}-\ref{abatwist}, we have, for $j$ and $k$ nonzero:

\begin{lem}\label{unique}$H_a^j=H_b^k$ $\Leftrightarrow$ $a=b$ and $j=k$.\end{lem}

\begin{lem}\label{disjoint}$[H_a^j,H_b^k]=1 \Leftrightarrow \gin(a,b) = 0$.\end{lem}

\begin{lem}\label{aba}Distinct half-twists $H_a$ and $H_b$ satisfy $H_a^jH_b^kH_a^j = H_b^kH_a^jH_b^k$ if and only if $a$ and $b$ are adjacent and $j=k=\pm 1$.\end{lem}

For this lemma, it suffices to note that the condition of adjacency of arcs in $\dpu$ is equivalent to the lifts of the arcs having intersection number 1.

\p{Remark.} Lemmas~\ref{word} through~\ref{aba} hold not only for $\bmodz$, but also in the contexts of $\bn$ and $\mods$ for $S$ a sphere with $n$ punctures:

{\em Braid groups.} If $n \geq 4$, two half-twists satisfy a relation in $\bn$ if and only if they satisfy the same relation in $\bmodz$.

{\em Punctured spheres.} If $W(H_a,H_b)$ is a word in the half-twists about $a$ and $b$ which represents the identity in $\mods$, then we can choose a puncture $p$ which is not an end of $a$ or $b$ (if $n \geq 5$).  Then, $W(H_a,H_b)$ is also the identity in the subgroup of $\mods$ consisting of maps which fix $p$.  But this subgroup is isomorphic to $B_{n-1}/\cent$.

\p{Remark.} Lemma~\ref{aba} is an algebraic characterization of geometric intersection number 2 (for 2-curves).  Another approach is to apply a theorem of Hamidi-Tehrani and the second author which says, in the case of nondisjoint curves $a$ and $b$ in $\dn$, that $T_a^jT_b^k$ is equal to a multitwist for $j,k \neq 0$ if and only if $\gin(a,b)=2$ and $j=k=\pm 1$ \cite{hht} \cite{dm}.

%%%
%%%
%%%

\section{$\bmodz$ is co-Hopfian}\label{bmodztobmodz}

In this section we prove Main Theorem~\ref{main3}, namely, that $\bmodz$ is co-Hopfian.  Throughout, we assume $n \geq 4$.

\subsection{Almost half-twist preserving}\label{almosttp}

Let $\inj$ be an injective endomorphism of $\bmodz$.  Our first goal is to show that $\inj$ is almost half-twist preserving; that is, for each half-twist, some power of it is mapped to a power of a half-twist.

The following result of Ivanov and McCarthy will be used:

\begin{thm}[Ivanov]\label{thm 5.9 I-M}Let $\Ga < \mods$ be a finite index subgroup consisting of pure elements.  If $f \in \Ga$ has canonical reduction system $\coll$, then \[ Z(C_\Ga(f)) \cong \z^{c+p} \] where $c$ is the number of curves in $\coll$ and $p$ is the number of pseudo-Anosov components of $f_\coll$. \end{thm}

This theorem is implicit in lecture notes of Ivanov \cite{ni}.  For a more recent exposition, refer to the paper of Ivanov and McCarthy \cite{im}.

\p{Setup.} For the remainder of the section, $\inj: \bmodz \to \bmodz$
is an injective homomorphism.  We also define: \[ \Ga' = \pbn/\cent \quad \text{ and } \quad \Ga = \inj^{-1}(\Ga') \cap \pbn/\cent \] Both of these subgroups consist entirely of pure elements of $\bmodz$ by Theorem~\ref{purity}.  Let $k=[\bmodz:\Ga]!$; note then $g^k \in \Ga$ for all $g \in \bmodz$.

\begin{prop}\label{thmatp}$\inj$ is almost half-twist preserving.\end{prop}

\bpf

Let $a$ be a 2-curve in $\dn$.  Then $f = H_a^k$ is an element of $\Ga$, and belongs to a maximal rank free abelian subgroup of $\bmodz$.

By Lemma~\ref{lem 10.2 I-M} and Theorem~\ref{thm 5.9 I-M},
$Z(C_{\Ga'}(\inj(f)))$ has rank at most $1$.  According to
Theorem~\ref{thm 5.9 I-M}, a canonical reduction system for $\inj(f)$
has $c$ circles and $p$ pseudo-Anosov components, where $c + p \leq
1$.  If $p=1$, then $\inj(f)$ is pseudo-Anosov, which contradicts the
fact that the centralizer of $f$ (and hence of $\inj(f)$) contains a
free abelian group of rank 2.  We can't have $c=p=0$, for then $\inj(f)$ is a finite order pure mapping class, and is hence the identity.

Thus, $c=1$ and $p=0$.  So there is nontrivial simple closed curve
$a'$ in $\dn$ such that $\inj(f) = T_{a'}^{k'}$ for some $k'$.  We now
show that $a'$ is a 2-curve.  Consider a maximal collection of
disjoint 2-curves in $\dn$, $\set{a=a_1, \dots, a_{\lfloor n/2
    \rfloor}}$.  The twists $\set{H_{a_i}}$ define a basis for a free
abelian group of rank $\lfloor n/2 \rfloor$, all of whose generators
are conjugate in $\bmodz$.  As $\inj$ is an injective homomorphism, $\set{\inj(H_{a_i}^k)}$ is a set of $\lfloor n/2 \rfloor$ twists about disjoint curves surrounding the same number of punctures.  Thus, all the curves are 2-curves.

\epf

\p{Action on curves.} By the above proposition and Lemma~\ref{unique}, $\inj$  has a well-defined action $\injs$ on 2-curves defined by: \[ \inj(H_a^k) = T_{\injs(a)}^{k'} \]
%Further, $\injs$ takes 2-curves to 2-curves.

%%%
%%%
%%%

\subsection{Half-twist preserving}\label{tp}

We now prove that any injective endomorphism of $\bmodz$ must be half-twist preserving, that is $\inj(H_a)=H_{a'}^{\pm 1}$ for any half-twist $H_a$.

\begin{prop}\label{htp}$\inj$ is half-twist preserving.\end{prop}

\bpf

Let $a$ be a 2-curve in $\dn$, and let $\injs(a)=a'$.  By Proposition~\ref{thmatp}, we have $\inj(H_a^k) = T_{a'}^{k'}$.  Since $[H_a,H_a^k]=1$ it follows that $[\inj(H_a),\inj(H_a^k)] =
[\inj(H_a),T_{a'}^{k'}] = 1$, and so $\inj(H_a)(a')=a'$ (Fact~\ref{fix}).

Let $S_1$ and $S_2$ be the surfaces obtained by cutting $\dn$ along the 2-curve $a'$.  If $S_1$ is the twice-punctured disk (the interior of $a'$), then $S_2$ is an annulus with $n-2$ punctures.

Since $\inj(H_a)$ fixes $a'$, there are well-defined reductions $f_1$ and $f_2$ of $\inj(H_{a})$ to $S_1$ and $S_2$.  These reductions must be finite order mapping classes since $\inj(H_a)^k = T_{a'}^{k'}$.

Since $S_1$ is a twice-punctured disk, $f_1$ must be a power of a half-twist.  To show that $f_2$ is the identity, we consider a 2-curve $b$ disjoint from $a$.  We have that $\injs(b)=b'$ is a 2-curve on $S_2$.  Further, by commutativity, $\inj(H_a)$ (and hence $f_2$) fixes $b'$.  If $S_2'$ is the complement of the interior of $b'$ in $S_2$, then we see that $f_2$ restricted to $S_2'$ is the identity (apply Theorem~\ref{eilenberg}).  It follows that $f_2$ is the identity.  Thus, we have: $\inj(H_a) = H_{a'}^m$.

% We observe that Lemma~\ref{unique} implies that $m = 2k'/k$.

Let $b$ be a 2-curve which is adjacent to $a$.  Since $H_b$ is conjugate to $H_a$, it follows that $\inj(H_b)=H_{b'}^{m}$, and then: \[ H_{a'}^{m}H_{b'}^{m} H_{a'}^{m} = H_{b'}^{m} H_{a'}^{m}H_{b'}^{m} \] By Lemma~\ref{aba}, $m=\pm 1$, and we are done.

\epf

%%%
%%%
%%%

\subsection{Homeomorphism}\label{homeo}

In this section we complete the last step in the proof of Main Theorem~\ref{main3}:

\begin{prop}$\inj$ is induced by a homeomorphism.  In particular, $\inj$ is an automorphism. \end{prop}

\bpf

Let $\set{a_1,\dots,a_{n-1}}$ be a sequence of arcs in $\dn$ with the property that $a_i$ is adjacent to $a_{i+1}$ for $1 \leq i \leq n-2$ and the arcs are disjoint and do not share ends otherwise.

By Proposition~\ref{htp} and Lemmas~\ref{disjoint}~and~\ref{aba}, the arcs $\set{\injs(a_i)}$ have the same properties.  Since $\dn - \set{a_1,\dots,a_{n-1}}$ is a punctured disk, there are exactly two mapping classes, say $h_+$ (orientation preserving) and $h_-$ (orientation reversing), whose actions on the $\set{a_i}$ agree with that of $\injs$.  Choose $h$ to be $h_+$ if $\inj(H_{a_1})$ is a positive half-twist, and $h=h_-$ otherwise.  Since half twists about the $a_i$ generate $\bmodz$, it follows that $\inj$ is induced by $h$.

\epf

%%%
%%%
%%%

\section{Proof of Main Theorem~\ref{main2}}\label{bntobnp}

Let $\inj : \bn \to \bnp$ be an injective homomorphism, and assume
$n \geq 7$ throughout this section.

\subsection{Almost half-twist preserving}\label{ahtp2} The following changes are made to Section~\ref{almosttp}: $\Ga'$ is now $\pbnp$, $\Ga = \pbn \cap \inj^{-1}(\Ga')$, and $k=[\bn:\Ga]!$.  Theorem~\ref{thm 5.9 I-M} is changed as follows:  all ranks of centers of centralizers increase by 1 since $\pbn$ has infinite cyclic center (replace $c+p$ with $c+p+1$).

Now, if $f \in \Ga$ is a power of a Dehn twist, then $Z(C_{\Ga'}(\inj(f)) \cong \z^{c+p+1}$, where $1 \leq c+p+1 \leq 3$ by Lemma~\ref{lem 10.2 I-M}, Theorem~\ref{thm 5.9 I-M}, and the fact $\bnp$ has an infinite cyclic center.  So in the current situation, there are more possibilities for $c$ and $p$.

For the following lemma, we say that pure mapping classes have {\em
 overlapping pseudo-Anosov components} if their reductions have
 pseudo-Anosov components which are distinct and nondisjoint.  By work
 of Ivanov \cite{nvi}, we have:

\begin{lem}\label{overlap}Maps with overlapping pseudo-Anosov components do not commute.\end{lem}

We say that multitwists {\em overlap} if any of the curves intersect.

\begin{lem}\label{overlap2}Overlapping multitwists do not commute.\end{lem}

The previous lemma follows from Lemma~\ref{commutecrs}, the definition of canonical reduction system, and the fact that the canonical reduction system for a multitwist is the set of curves in the multitwist (the last fact is due to Birman, Lubotzky, and McCarthy \cite{blm}).

\begin{prop}\label{ahtp}If $f=H_a^k \in \Ga$, then $\inj(f)$ is the product of a multitwist (about at most two curves) and a central element of $\bnp$.\end{prop}

\bpf

As in Section~\ref{almosttp}, we have: \[ 1 \leq c + p + 1 \leq 3 \] where $c$ is the number of components in the canonical reduction system for $\inj(f)$, and $p$ is the number of its pseudo-Anosov components.  Thus, we have the following possibilities for $(c,p)$: \[ (0,2) \quad (0,1) \quad (0,0) \quad (1,1) \quad (2,0) \quad (1,0) \] The first possibility is absurd.  The second and third possibilities are ruled out for the same reasons as before.  The goal is to show that only the last two possibilities occur, so it remains to rule out the fourth.

Assume $c=p=1$.  Pick a maximal collection of disjoint 2-curves $\set{a=a_1, \dots, a_{\lfloor n/2 \rfloor}}$.  The $\inj(H_{a_i}^k)$ are all conjugate, so they each have one curve $a_i'$ in their canonical reduction system, and one pseudo-Anosov component.  Since they all commute, the $a_i'$ are disjoint by Lemma~\ref{commutecrs}.  So the $a_i'$ must all be 2-curves (using conjugacy).
But the pseudo-Anosov pieces cannot be in the interior of the $a_i'$ (there are no pseudo-Anosov maps), and they cannot be in the exterior of the $a_i'$, because then they all overlap, violating Lemma~\ref{overlap}.

\epf

\p{Action on curves.} By Proposition~\ref{ahtp}, there is a function: \[ \injs: \set{\text{2-curves}} \to \set{\text{multicurves}} \] By {\em multicurve}, we mean a collection of mutually disjoint nontrivial curves.

We note that since $\inj$ is an injective homomorphism,
Lemma~\ref{overlap2} implies that $\injs$ preserves disjointness and
non-disjointness of collections of curves.  We are now ready to show
that $\inj$ is almost half-twist preserving in the sense of Section~\ref{almosttp}.

\begin{prop}\label{ahtpprop2}If $f=H_a^k \in \Ga$, then $\inj(f) = H_{a'}^{k'} H_A^s \gar^t$ for some arc $a'$ and curve $A$, where $H_A^s \gar^t$ is independent of $a$.\end{prop}

\bpf

If $\injs(a)$ is a single curve $a'$, then $\inj(f) = H_{a'}^{k'}z^t$, for some $t$.  By conjugacy, $A$, $s$, and $t$ are independent of $a$ ($s = 0$).  Now assume $\injs(a)=\set{a',A}$.  By conjugacy then, every half-twist then maps to a multitwist about two curves (modulo central elements).

First, by combining Theorem~\ref{thm 5.9 I-M} and Lemma~\ref{lem 10.2 I-M}, we see that if $\gin(a,b)=0$ for 2-curves $a$ and $b$, then $\injs(a)$ and $\injs(b)$ must share a curve, say $\injs(b)=\set{b',A}$.

Now we show that if $a$, $b$, and $c$ are mutually disjoint 2-curves, then the multicurves $\injs(a)$, $\injs(b)$, and $\injs(c)$ share a common curve (namely, $A$).  Suppose not.  Then $\injs(c) = \set{a',b'}$.  This configuration contradicts the fact that there are elements which commute with $\inj(H_a^k)$ and $\inj(H_b^k)$ but not $\inj(H_c^k)$ (apply Fact~\ref{fix}).

Now we want to show that if $x$ is any 2-curve, we have $\injs(x)=\set{x',A}$.  We will see that this is implied by the connectedness of the following graph, which we call the {\em arc triple complex} $\atdn$: \bl \item[Vertices.] Triples of disjoint 2-curves (thought of as arcs)  \item[Edges.] Two triples which have two curves in common \el As above, associated to each vertex $v$ of $\atdn$ is a unique curve $A_v$ (the one common to the images under $\injs$ of the three arcs).  If $v$ is connected to $w$ by an edge, we see that $A_w=A_v$.  Since $\atdn$ is connected for $n \geq 7$ (see Section~\ref{ad}), there is a unique curve $A$ so that $A \in \injs(x)$ for all $x$.

As in Section~\ref{almosttp}, $a'$ must be a 2-curve: any maximal collection $\set{a=a_1,\dots,a_{\lfloor n/2 \rfloor}}$ of disjoint 2-curves must go to a set of $\lfloor n/2 \rfloor$  2-curves (plus possibly $A$), since the corresponding abelian subgroup with conjugate generators is preserved.

\epf

\subsection{Half-twist preserving} We now show that $\inj$ is half-twist preserving in the sense of Section~\ref{tp}.

\begin{prop}\label{ngeq7prop}If $H_a \in \bn$ is a half-twist, then
  $\inj(H_a) = H_{a'}^{\pm 1} f_a H_A^s \gar^t$ for a 2-curve $a' \in
  \injs(a)$, where $H_A^sz^t$ is independent of $a$, and $f_a$ is
  supported on the component of $\dnp-A$ which does not contain $a'$.\end{prop}

\bpf

Let $H_a$ be a half-twist.  If $\injs(a)$ is a single curve (in this case $s = 0$), then the argument in Section~\ref{tp} applies.  Thus, we assume $\injs(a)=\set{a',A}$, where $a'$ and $A$ are as in Proposition~\ref{ahtpprop2}.

As in Section~\ref{tp}, we see that $\inj(H_a)$ fixes $\set{a',A}$ by commutativity (Fact~\ref{fix}).  If $b$ is any 2-curve in $\dn$ disjoint from $a$, we have $\injs(b)=\set{b',A}$, and it follows from the commutativity of $H_a$ and $H_b$ that $\inj(H_a)$ also fixes $\set{b',A}$, and hence fixes $a'$ and $A$ individually.

Thus, we can consider the reductions of $\inj(H_a)$ to the surfaces
obtained by cutting $\dnp$ along $a'$ and $A$.  Let $S_1$ be the
interior of $a'$, let $S_2$ be the component containing $b'$, and let
$S_3$ be the last component.  By conjugacy of $H_a$ and $H_b$, we have
that $a'$ is a boundary component of $S_2$.  Again, these reductions
are finite order mapping classes.

As before, the mapping class $f_1$ of $S_1$ must be a power of a
half-twist.  By commutativity, the mapping class $f_2$ of $S_2$ fixes
$b'$.  Since it also fixes $a'$ and $A$, then as in Section~\ref{tp}
it is the identity.

Hence $\inj(H_a)=H_{a'}^mf_aH_A^s\gar^t$, where $f_a$ is some mapping class supported on $S_3$.

For an arc $b$ adjacent to the arc $a$, we have that $H_a$ and $H_b$ satisfy the braid relation, and $\inj(H_b)=H_{b'}^mf_bH_A^sz^t$.  By considering the reductions of $\inj(H_a)$ and $\inj(H_b)$ to the surface $S_1 \cup S_2$, we have: \[ H_{a'}^m H_{b'}^m H_{a'}^m = H_{b'}^m H_{a'}^m H_{b'}^m \]
as elements of $\mcg(S_1 \cup S_2)$.  Lemma~\ref{aba} implies that $m = \pm 1$.

\epf

\subsection{Embedding}

The following is the last step of the proof of Main Theorem~\ref{main2}:

\begin{prop}$\inj$ is induced by an embedding $h: \dn \to \dnp$.\end{prop}

\bpf

Consider a chain of arcs $\set{a_1,\dots,a_{n-1}}$ connecting successive punctures in $\dn$.  By Proposition~\ref{ngeq7prop} and Lemmas~\ref{disjoint} and~\ref{aba}, we have that the $\set{a_i'}$ (where $\injs(a_i) = \set{a_i'}$ or $\set{a_i',A}$) is a similar chain of arcs connecting $n$ of the punctures in $\dnp$.  Choose an embedding $h: \dn \to \dnp$ taking the first chain to the second chain, and taking $\ddn$ to the boundary of a regular neighborhood of the second chain (as in Section~\ref{homeo}, there are two choices).  First note that if $A \in \injs(a_i)$, then $A = h(\ddn)$, for otherwise it would intersect one of the $a_i'$ or be isotopic to $\partial \dnp$ (it follows that each $f_{a_i}$ is trivial).

We see that $h$ induces $\inj$ since $\inj(H_a) = H_{h(a)}^{\ep(h)}T_A^s\gar^t$ and the $\set{H_{a_i}}$ generate $\bn$.  Note that $H_A^s$ must actually be a power of a full Dehn twist since it commutes with the images of the generators.

\epf

\p{Special cases.}  Using ideas of Chris Leininger, we have complete
proofs of Main Theorem~\ref{main2} for $n \in \set{4,5,6}$.  For $n=5$
and $n=6$, we use the connectedness of $A_2(D_n)$; the key is to
consider elements of $\bn$ which ``realize'' the edges of the complex.
The case of $n=4$ requires more work, since $A_2(D_4)$ is not
connected (it has infinitely many components).

%%%
%%%
%%%

\section{Mapping class groups of spheres}\label{mcgpf}

We now explain why the proofs of the main theorems also apply to the case of mapping class groups of spheres (as per Theorem~\ref{mcgthm}).

Let $S$ be a sphere with $n \geq 5$ punctures.  We first note that $B_{n-1}/\cent$ is a finite index subgroup of $\mods$ (the subgroup fixing one of the punctures).  Let $S'$ be either $S$ or a sphere with $n+1$ punctures.

To see that any injective homomorphism $\inj: \mods \to \mcg(S')$ is almost half-twist preserving, we study finite index subgroups $\Ga < \mods$ and $\Ga' < \mcg(S')$, which are the same groups as in the braid case.

If $S=S'$, the proofs that $\inj$ is half-twist preserving and that $\inj$ induced by a homeomorphism are the same as for the braid case.

For $S \neq S'$, we use the obvious spherical version of the arc k-tuple complex, and, as in the braid case, we find that $\inj(H_a)=H_{a'}^mH_A^s$ for any half-twist $H_a$, where $a'$ is a 2-curve and $H_A^s$ is independent of $a$.  If $b$ is adjacent to $a$, then an application of Lemma~\ref{aba} yields that $m=\pm 1$, and if $\inj(H_b)=H_{b'}^mH_A^s$, that $a'$ is adjacent to $b'$.  Also, any chain $\set{a_i}$ of $n-1$ arcs in $S$ gets mapped to a chain $\set{a_i'}$ of $n-1$ arcs in $S'$.  Since $A$ is disjoint from this chain in $S'$, it follows that $A$ is trivial, that is, $\inj(H_a) = H_{a'}^{\pm 1}$.

The element $(H_{a_1}\cdots H_{a_{n-2}})^{n-1}$ is trivial in $\mods$.
But this gets mapped to $(H_{a_1'}\cdots H_{a_{n-2}'})^{n-1}$, which
is a nontrivial twist.  This is a contradiction, so there are no
injective homomorphisms $\mods \to \mcg(S')$.

%%%
%%%
%%%

\section{Arc k-tuple complex}\label{ad}

We define the following complex $\adnk$ for $k \geq 2$ and $n \geq 2k$:
\begin{list}{}{\setlength{\leftmargin}{1in}} \item[Vertices.]$k$-tuples of disjoint arcs in $\dn$ (all $2k$ ends distinct) \item[Edges.]Pairs of $k$-tuples sharing a $(k-1)$-tuple \end{list}

In Section~\ref{bntobnp}, we use the fact that the {\em arc triple complex} ($k = 3$) is connected for $n \geq 7$.
It is not hard to see that $\adnk$ is not connected for $n=2k$.

We think of $\dn$ as a disk in $\br^2$ with $n$ points removed along a line.  This allows us to speak of the {\em closest} puncture(s) to a given puncture.  A \emph{straight arc} is a linear arc (whose ends are necessarily neighboring punctures).

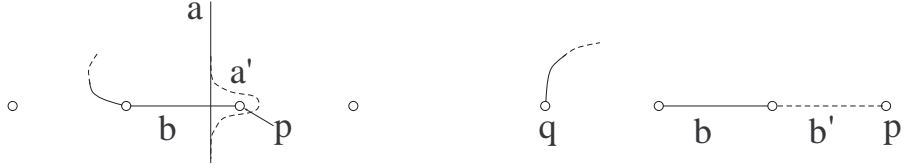
\pics{semihoriz}{Left: a basic move;  right: shuffling.}{3}

\p{Basic moves.} We will make frequent use of the following move which reduces the intersection between a vertex $A$ and an arc $b$.  Refer to the left diagram in Figure~\ref{semihorizpic}.  Suppose $p$ is an end of $b$, but not the end of any arc of $A$, and that $a$ is the arc of $A$ closest to $p$ along $b$.  The replacement of $a$ by $a'$, obtained by pushing $a$ off $b$, is an edge from $A$ to $A' = (A -\set{a}) \cup \set{a'}$ in $\adnk$.  We call this a {\em basic move} along $b$ at $p$.

\begin{thm}\label{thm: arcconnect} $\adnk$ is connected for $n > 2k$.\end{thm}

\bpf

Choose a base vertex $V$ of $\adnk$, all of whose arcs are straight, and let $A$ be any other vertex.  Assuming $A$ has at least one non-straight arc, we will show there is a path from $A$ to a vertex which has one more straight arc than $A$.  By induction, this implies that $A$ is connected to a vertex consisting only of straight arcs, and it is clear that any such arc is connected to $V$.

Because $n > 2k$, we can choose a puncture $p$ which is not an end of any arc in $A$ and which is closest to the set of ends of the non-straight arcs of $A$; let $q$ be one of these ends closest to $p$.  Potentially, there are straight arcs of $A$ between $p$ and $q$.  We shuffle these straight arcs as follows (refer to the right side of Figure~\ref{semihorizpic}):

Suppose that $b \in A$ is the straight arc closest to $p$.  Let $b'$ be the straight arc joining $p$ and the end of $b$ which is closest to $p$.  Perform basic moves along $b'$ at $p$ until the new arcs do not intersect $b'$.  Now replace $b$ with $b'$.  The end of $b$ furthest from $p$ is not the end of any arc in the new vertex and is closer to $q$ than $p$.  Repeating this process if necessary, we can assume that $p$ neighbors the end $q$ of a non-straight arc, say $a$, of $A$.

Let $c$ be the straight arc joining $p$ and $q$.  Perform basic moves along $c$ at $p$ until the new arcs do not intersect $c$.  Replacing $a$ with $c$ completes the inductive step.

\epf

%%%
%%%
%%%

\section{Questions}\label{q}

The results of this paper suggest a variety of directions for further
study.  Some answers have been obtained after the circulation of the
first version of this paper (see footnotes).

\begin{q}Is every injective homomorphism of $\bn$ into $B_m$ geometric?\end{q}

In general, injective homomorphisms are more complicated than in our
main theorems.  For example, consider the following two injective homomorphisms of $\bn$ into $B_{2n}$: \bl \item cabling: half twists are sent to elements switching {\em pairs} of punctures (see Figure~\ref{ghtpic}) \item doubling: half twists are sent to products of two half-twists (one in each ``half'' of $D_{2n}$) \el

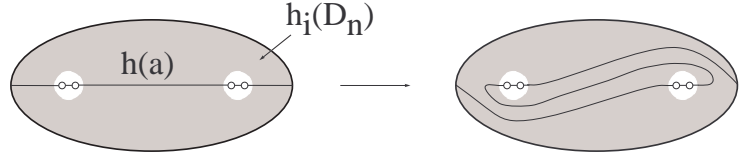
\pics{ght}{Generalized half-twist, $k=2$.}{3}

We now give a definition of {\em geometric injection}:

Let $\set{h_i}$ be a finite collection of embeddings $\dn \to D_m$ with: \bl \item $h_i(\dn)$ are mutually disjoint (possibly nested) \item for a given $i$, any circle with a single puncture in its interior is sent by $h_i$ to a circle with $k$ punctures in its interior, where $k$ does not depend on the circle \el Then, a {\em strictly geometric} injection is one defined as follows on generators: \[ \inj(H_a) = H_{h(a)} \] Here, $a$ is
any 2-curve and $H_{h(a)} = \prod H_{h_i(a)}$.  Note that in general $h_i(a)$ is not a  2-curve, but a $2k$-curve ($k$ as above).  In this case, $H_{h_i(a)}$ is a generalized half-twist, as in Figure~\ref{ghtpic}.

Note that both of the above examples of injective homomorphisms (cabling and doubling) are strictly geometric in this sense.

Actually, a strictly geometric injection can be augmented to an
injective homomorphism (still deserving of the name geometric) by adding appropriately defined ``constant terms'', as in our main theorems.  An example of a constant term is a product of twists about disjoint curves which are themselves disjoint from the $h_i(\dn)$.

\p{Braiding braids.} Nested embeddings of disks give rise to even more
interesting strictly geometric injections.  For instance, an injective
homomorphism of $\bn$ into $B_{n^2}$
is given by including $n$ copies of $\bn$ into $B_{n^2}$, and also braiding the $n$ copies.

For an appropriate definition of geometric injection, we have:

\begin{q}[Farb--Margalit]Is every injection of mapping class groups geometric?\end{q}

The co-Hopfian property has not been established for certain related groups:

\begin{q}Which Artin groups modulo their centers
  co-Hopfian?\footnote{For partial results, see \cite{bellm}.}\end{q}

\begin{q}Are surface braid groups co-Hopfian?\end{q}

Irmak, Ivanov, and McCarthy have recently shown that the automorphism group of a higher-genus surface braid group is the extended mapping class group of the corresponding punctured surface \cite{iim}.

Superinjective maps of curve complexes were introduced by Irmak in
order to study injective homomorphisms of finite-index subgroups of mapping class groups (for higher genus surfaces) \cite{ei}.

\begin{q}\label{si}Is every superinjective map of the complex of
  curves of a punctured sphere or punctured torus induced by a homeomorphism of the
  surface?\footnote{This is completely answered in \cite{bm} \cite{bellm}
  \cite{kjs}.} \end{q}

An affirmative answer to Question~\ref{si} would extend the results of this paper to finite index subgroups, in particular pure braid groups.  Along the same lines, we have:

\begin{q}What is the abstract commensurator of
  $\bn$?\footnote{Complete answer given in \cite{lm}.}\end{q}

%%%
%%%
%%%

\bibliographystyle{plain}

\bibliography{cohopf}

\end{document}